\documentclass[12pt]{article}
\usepackage{enumerate, amsmath, amsthm, amsfonts, amssymb, xy}
\input xy
\xyoption{all}
\usepackage[margin=1in]{geometry}
\usepackage[dvips]{graphicx,color}
\usepackage{graphicx}
\usepackage{overpic}
\usepackage{multimedia}
\usepackage{stmaryrd}
\usepackage{amsfonts}
\usepackage{pifont}
\usepackage{amssymb,amsmath,amsfonts,mathrsfs,amscd,pict2e}
\usepackage{flushend, cuted}
\usepackage{multicol}
\usepackage{color}
\usepackage{longtable,colordvi}
\newtheorem{thm}{Theorem}[section]
\newtheorem{prop}[thm]{Proposition}
\newtheorem{defi}[thm]{Definition}
\newtheorem{lem}[thm]{Lemma}
\newtheorem{cor}[thm]{Corollary}
\newtheorem{exa}[thm]{Example}

\def\qed{\nopagebreak\hfill{\rule{4pt}{7pt}}\medbreak}
\def\pf{\noindent{\it Proof.} }

\makeatletter \@addtoreset{equation}{section} \makeatother

\def\Cc{{\mathcal C}}
\def\Oc{{\mathcal O}}

\def\bold{\bf}

\def\eb{{\bold e}}

\def\0b{{\bold 0}}

\def\Cc{{\mathcal C}}

\usepackage{bm}
\bmdefine{\Bzero}{0}
\bmdefine{\Bone}{1}
\def\Bone{{\bf 1}}
\def\RR{{\mathbb R}}

\begin{document}
\begin{center}
{\Large\bf Order-Chain Polytopes}
\end{center}
\renewcommand{\thefootnote}{\fnsymbol{footnote}}
\begin{center}
Takayuki Hibi$^{a}$, Nan Li$^{b}$, Teresa Xue Shan Li\footnote{Corresponding author.}$^{c}$,  Li Li Mu$^{d}$ and Akiyoshi Tsuchiya$^{a}$\\[6pt]
$^a$Department of Pure and Applied Mathematics, Graduate School of Information Science and Technology, Osaka University, Toyonaka,
Osaka 560-0043, Japan\\[6pt]
$^b$Department of Mathematics, Massachusetts Institute of Technology\\
 Cambridge, MA 02139, USA\\[6pt]
$^{c}$School of Mathematics and Statistics, Southwest University\\
  Chongqing 400715, PR China\\[6pt]
$^{d}$ School of Mathematical Sciences,
 Dalian University of Technology\\
  Dalian 116024, PR China\footnote[0]{E-mail addresses: hibi@math.sci.osaka-u.ac.jp (T. Hibi), amenda860111@gmail.com (N. Li),\\
pmgb@swu.edu.cn (T.X.S.Li), lly-mu@hotmail.com (L.L. Mu),\\ a-tsuchiya@cr.math.sci.osaka-u.ac.jp(A. Tsuchiya)} \\[6pt]

\end{center}

\begin{abstract}
Given two families $X$ and $Y$ of integral polytopes with nice combinatorial and algebraic properties, a natural way to generate new class of polytopes is to take the intersection $\mathcal{P}=\mathcal{P}_1\cap\mathcal{P}_2$, where $\mathcal{P}_1\in X$, $\mathcal{P}_2\in Y$.  Two basic questions then arise:  1) when $\mathcal{P}$ is integral and 2) whether $\mathcal{P}$ inherits the ``old type" from $\mathcal{P}_1, \mathcal{P}_2$ or  has a ``new type", that is, whether $\mathcal{P}$ is unimodularly equivalent to some polytope in $X\cup Y$ or not. In this paper, we focus on the families of order polytopes and chain polytopes and create a new class of polytopes following the above framework, which are named order-chain polytopes. In the study on their volumes, we discover a natural relation with Ehrenborg and Mahajan's results on maximizing descent statistics.
\end{abstract}

\noindent\textbf{Keywords:} poset, order-chain polytope, unimodular equivalence

\section{Introduction}
This paper was motivated by the following two questions on intersecting two integral polytopes: $\mathcal{P}_1\in X$, $\mathcal{P}_2\in Y$, where $X$ and $Y$ are two families of polytopes with nice combinatorial and algebraic properties:  1) when the intersection $\mathcal{P}=\mathcal{P}_1\cap\mathcal{P}_2$ is integral and 2) whether $\mathcal{P}$ inherits the ``old type" from $\mathcal{P}_1, \mathcal{P}_2$ or  has a ``new type", that is, whether $\mathcal{P}$ is unimodularly equivalent to some polytope in $X\cup Y$ or not. In this paper,  we mainly consider the families of order polytopes and chain polytopes. Instead of taking the intersection of an arbitrary $d$-dimensional order polytope and an arbitrary $d$-dimensional chain polytope, we will consider the intersection of an order polytope $\mathcal{O}(P')$ and a chain polytope $\mathcal{C}(P'')$, both of which arise from subposets $P', P''$ of a given poset. This leads us to the notion of order-chain polytope, which generalizes both order polytope and chain polytope.

The order polytope ${\mathcal O}(P)$ as well as the chain polytope
${\mathcal C}(P)$ arising from a finite partially ordered set $P$
has been studied by many authors from viewpoints of both combinatorics
and commutative algebra.  Especially, in Stanley \cite{Stanley1986},
the combinatorial structure of order polytopes and chain polytopes
is explicitly discussed.  Furthermore, in \cite{Hibi-Li2012},
the natural question when the order polytope ${\mathcal O}(P)$ and
the chain polytope ${\mathcal C}(P)$ are unimodularly equivalent is
solved completely.  It follows from \cite{Hibi} and \cite{ASL} that
the toric ring (\cite[p.~37]{JSTgb}) of ${\mathcal O}(P)$ and
that of ${\mathcal C}(P)$ are algebras with straightening laws
(\cite[p.~124]{HibiRedBook}) on finite distributive lattices.
Thus in particular the toric ideal (\cite[p.~35]{JSTgb}) of each of
${\mathcal O}(P)$ and ${\mathcal C}(P)$ possesses a squarefree quadratic
initial ideal (\cite[p.~10]{JSTgb}) and possesses a regular unimodular
triangulation (\cite[p.~254]{JSTgb}) arising from a flag complex.
(Recall that a flag complex is a simplicial complex any of its nonface
is an edge.)  Furthermore, toric rings of order polytopes naturally
appear in algebraic geometry (e.g., \cite{BL})
and in representation theory (e.g., \cite{Wang}).

We begin by introducing some basic notation and terminology. Given a convex polytope $\mathcal{P} \subset \mathbb{R}^{d}$, we write
$\mathcal{V}(\mathcal{P})$ for the set of vertices of $\mathcal{P}$ and
$\mathcal{E}(\mathcal{P})$ for the set of edges of $\mathcal{P}$.
A {\em facet hyperplane} of $\mathcal{P} \subset \mathbb{R}^{d}$ is defined to be a hyperplane
of $\mathbb{R}^{d}$ which contains a facet of $\mathcal{P}$.
If
\[
H=\{ \, (x_1,x_2,\ldots, x_d) \in \mathbb{R}^{d} \, : \, a_1x_1+a_2x_2+\cdots+a_dx_d-b=0 \, \},
\]
where each $a_{i}$ and $b$ belong to $\mathbb{R}$,
is a hyperplane of $\mathbb{R}^{d}$ and $v=(y_1,y_2,\ldots,y_d)\in \mathbb{R}^{d}$,
then we set
\[
H(v)=a_1y_1+a_2y_2+\ldots+a_dy_d-b.
\]

Let $(P,\preccurlyeq)$ be a finite partially ordered set ({\em poset}, for short)
on $[d]=\{1,\ldots,d\}$.
For each subset $S\subseteq P$, we define
$\rho(S)=\sum_{i\in S}\eb_i$, where $\eb_1,\ldots, \eb_d$ are the canonical unit
coordinate vectors of $\mathbb{R}^{d}$.
In particular $\rho(\emptyset)=(0,0,\ldots,0)$, the origin of $\mathbb{R}^{d}$.
A subset $I$ of $P$ is an {\em order ideal} of $P$
if $i\in I$, $j \in [d]$ together with $j\preccurlyeq i$ in $P$ imply $j\in I$.
An {\em antichain} of $P$ is a subset $A$ of $P$ such that any two elements
in $A$ are incomparable.  We say that $j$ {\em covers} $i$
if $i\prec j$ and there is no $k\in P$ such that $i\prec k\prec j$.
A chain $j_1\prec j_2\prec\cdots\prec j_s$ is {\em saturated}
if $j_q$ covers $j_{q-1}$ for $1<q\leq s$.
A poset can be represented with its Hasse diagram, in which each cover relation $i\prec j$ corresponds to an edge denoted by $e=\{i,j\}$.

In \cite{Stanley1986}, Stanley introduced two convex polytopes arising from
a finite poset, the order polytope and the chain polytope.
Following \cite{Hibi-Li2012}, we employ slightly different definitions.
Given a finite poset $(P,\preccurlyeq)$ on $[d]$,
the {\em order polytope} $\mathcal{O}(P)$ is defined to be the convex polytope
consisting of those $(x_1,\ldots,x_d)\in\mathbb{R}^{d}$ such that
\begin{itemize}
\item[(1)] $0\leq x_i\leq 1$ for $1\leq i\leq d$;\\
\item[(2)] $x_i\geq x_j$ if $i\preccurlyeq j$ in $P$.
\end{itemize}
The {\it chain polytope} $\mathcal{C}(P)$ of $P$ is defined to be the convex polytope
consisting of those $(x_1,\ldots,x_d)\in\mathbb{R}^{d}$ such that
\begin{itemize}
\item[(1)] $x_i\geq 0$ for $1\leq i\leq d$;\\
\item[(2)] $x_{i_1}+\cdots+x_{i_k}\leq 1$
for every maximal chain $i_1\prec\cdots\prec i_k$ of $P$.
\end{itemize}

Let $P$ be a finite poset and $E(P)$ the set of edges of its Hasse diagram.
In the present paper, an {\em edge partition} of $P$ is a map
\[
\ell: E(P)\longrightarrow \{o,c\}.
\]
Equivalently, an edge partition of $P$ is an ordered pair
\[
(oE(P), cE(P))
\]
of subsets of $E(P)$ such that $oE(P)\cup cE(P)=E(P)$ and $oE(P)\cap cE(P)=\emptyset$.
An edge partition $\ell$ is called {\em nontrivial} if
$oE(P) \neq \emptyset$ and $cE(P) \neq \emptyset$.

Suppose that $(P,\preccurlyeq)$ is a poset on $[d]$ with an edge partition
$\ell=(oE(P),cE(P))$.
Let $P_\ell^{'}$ and $P_\ell^{''}$ denote the $d$-element subposets of $P$ with edge sets $oE(P)$ and $cE(P)$ respectively. The {\em order-chain polytope} $\mathcal{OC}_{\ell}(P)$
with respect to the edge partition $\ell$ of $P$
is defined to be the convex polytope
\[
\mathcal{O}(P_\ell^{'})\cap \mathcal{C}(P_\ell^{''})
\]
in $\RR^{d}$.
Clearly the notion of order-chain polytope is a natural generalization of
both order polytope and chain polytope of a finite poset.

For example, let $P$ be the chain $1 \prec 2 \prec \cdots \prec 7$
with
\[
oE(P) = \{\{1,2\}, \{4,5\}, \{5,6\}\}, \, \, \,
cE(P) = \{\{2,3\}, \{3,4\}, \{6,7\}\}.
\]
Then $P_\ell^{'}$ is the disjoint union of the following four chains:
$$1 \prec 2,\ \ 3,\ \ 4 \prec 5 \prec 6,\ \ 7$$
and $P_\ell^{''}$ is the disjoint union of
$$1,\ \  2 \prec 3 \prec 4,\ \  5\ \ {\rm and}\ \ 6 \prec 7.$$
Hence the order-chain polytope $\mathcal{OC}_{\ell}(P)$ is the convex polytope
consisting of those $(x_1,\ldots,x_7)\in\mathbb{R}^{7}$ such that
\begin{itemize}
\item[(1)] $0\leq x_i\leq 1$ for $1\leq i\leq 7$;\\
\item[(2)] $x_1\geq x_2,\  x_{4}\geq x_{5} \geq x_{6}$;\\
\item[(3)] $x_{2}+x_{3}+x_{4}\leq 1,\  x_{6}+x_{7}\leq 1$.
\end{itemize}

 It should be noted that, for any poset $P$ on $[d]$ and any edge partition $\ell$ of $P$, the dimension of the order-chain polytope $\mathcal{O}\mathcal{C}_\ell(P)$ is equal to $d$. In fact, let $x=(1/d,\ldots,1/d) \in \mathbb{R}^d$, clearly, we have $x\in\mathcal{O}\mathcal{C}_\ell(P)$.  If $P'_\ell$ is an antichain, then $\mathcal{O}(P'_\ell)$ is the $d$-cube $[0,1]^d$. In this case,  $\mathcal{O}\mathcal{C}_\ell(P)$ is exactly the same as the chain polytope $\mathcal{C}(P)$ and so is $d$-dimensional.  If  $P'_\ell$ is not an antichain, then $P''_\ell$ is not a $d$-element chain. In this case, $x\in\partial\mathcal{O}(P'_\ell)$ and  $x\in\mathcal{C}(P''_\ell)\setminus\partial\mathcal{C}(P''_\ell)$, since no facet hyperplane of $\mathcal{C}(P''_\ell)$ contains $x$. In this case, we can find a ball $B_d(x)$ centered at $x$ such that $B_d(x)\subset \mathcal{C}(P''_\ell)\setminus\partial\mathcal{C}(P''_\ell)$. Keeping in mind that $x$ belongs to the boundary of $\mathcal{O}(P'_\ell)$, we deduce that $B_d(x)\cap (\mathcal{O}(P'_\ell)\setminus\partial\mathcal{O}(P'_\ell))\neq \emptyset$. It follows that $(\mathcal{O}(P'_\ell)\setminus\partial\mathcal{O}(P'_\ell))\cap (\mathcal{C}(P''_\ell)\setminus\partial\mathcal{C}(P''_\ell))\neq\emptyset$, as desired.

Recall that an integral convex polytope (a convex polytope is {\em integral} if all of its vertices have
integer coordinates) is called {\em compressed}
(\cite{rational}) if all of its ``pulling triangulations'' are unimodular.
Equivalently, a compressed polytope is an integral convex polytope
any of whose reverse lexicographic initial ideals is squarefree
(\cite{Stu}).
It follows from \cite[Theorem 1.1]{OhHicompressed} that all order polytopes
and all chain polytopes are compressed.  Hence the intersection
of an order polytope and a chain polytope is compressed if it is integral.
In particular every integral order-chain polytope is compressed.
It then follows that every integral order-chain polytope possesses
a unimodular triangulation and is normal (\cite{OhHinormal}).

Then one of the natural question, which we study in Section 2, is
when an order-chain polytope is integral.
We call an edge partition $\ell$ of a finite poset $P$ {\em integral} if
the order-chain polytope $\mathcal{OC}_{\ell}(P)$ is integral.
We show that every edge partition of a finite poset $P$
is integral if and only if $P$ is acyclic in Section 2. Here by an acyclic poset $P$ we mean that the Hasse diagram of $P$ is an acyclic graph.
Furthermore, we prove that every poset $P$ with $|E(P)|\geq 2$ possesses
at least one nontrivial integral edge partition.

In Section 3, we consider the problem when an integral order-chain polytope is unimodularly equivalent to either an order polytope or a chain polytope. This problem is related to the work \cite{Hibi-Li2012}, in which the authors characterize all finite posets $P$ such that $\mathcal{O}(P)$ and $\mathcal{C}(P)$ are unimodularly equivalent.  We show that if $P$ is either a disjoint union of chains or a zigzag poset, then the order-chain polytope
$\mathcal{OC}_{\ell}(P)$, with respect to each edge partition $\ell$ of $P$, is unimodularly equivalent to the chain polytope
of some poset (Theorems \ref{Chain} and \ref{Zigzag}).
On the other hand, for each positive integer $d\geq 6$,  we find a $d$-dimensional integral order-chain polytope which is not unimodularly equivalent to any chain polytope nor order polytope. This means that the notion of order-chain polytope is a nontrivial generalization of order polytope or chain polytope.

We conclude the present paper with an observation on the volume of
order-chain polytopes in Section 4.
A fundamental question is to find an edge partition $\ell$ of
a poset $P$ which maximizes the volume of $\mathcal{OC}_{\ell}(P)$.
In general, it seems to be very difficult to find a complete answer.
We discuss the case when $P$ is a chain on $[d]$, which involves Ehrenborg and Mahajan's problem (see \cite{Ehrenborg-Mahajan1998}) of maximizing the descent statistics over certain family of subsets.

\section{Integral order-chain polytopes}

In this section, we consider the problem when an order-chain polytope is integral. We shall prove that every edge partition of a poset $P$ is integral if and only if the poset $P$
is acyclic. We also prove that every poset $P$ with $|E(P)|\geq 2$ has at least one nontrivial integral edge partition.

\begin{thm}\label{THM1}
Let $P$ be a finite poset. Then every edge partition  of $P$ is integral if and only if $P$ is an acyclic poset.
\end{thm}
\pf Suppose that each edge partition $\ell$ of $P$ is integral. If the Hasse diagram of $P$ has a cycle $c$, then it is easy to find a non-integral edge partition. In fact, let $e=\{i,j\}$ be an arbitrary edge from $c$ and $\ell=(E(P)\setminus\{e\},\{e\})$. We now show that $\ell$ is not integral. To this end, let $I$ be the connected component of the Hasse diagram of $P_l^{'}$ which contains $i$ and $j$ and let $v=(v_1,v_2,\ldots,v_d)\in\mathbb{R}^d$ with
\[
v_k=\left\{
      \begin{array}{ll}
        \frac{1}{2}, & \hbox{if\ $k\in I$;} \\[5pt]
        0, & \hbox{otherwise.}
      \end{array}
    \right.
\]Then it is easy to see that
\[
 v=\bigcap_{\{p,q\}\in E(I)}H_{pq}\bigcap_{t\notin I}H_t\bigcap H_{ij},
 \]where
\begin{align*}
H_{pq}=\{(x_1,x_2,\ldots,x_d)\ |\  x_p=x_q\}\ \ {\rm for\ }e=\{p,q\}\in E(I) \\[5pt]
 H_t=\{(x_1,x_2,\ldots,x_d)\ |\ x_t=0\} \ \ {\rm for\ }t\notin I\\[5pt]
 H_{ij}=\{(x_1,x_2,\ldots,x_d)\ |\  x_i+x_j=1\}.
\end{align*}
are all facet hyperplanes of $\mathcal{OC}_{\ell}(P)$. So we deduce that $v$ is a
vertex of $\mathcal{OC}_{\ell}(P)$, a contradiction.

Conversely, suppose that $P$ is an acyclic poset on $[d]$  and $\ell$ is an edge partition of $P$. If $v=(a_1,a_2,\ldots,a_d)$ is a vertex of $\mathcal{OC}_\ell(P)$, then we can find $d$ independent facet hyperplanes of $\mathcal{OC}_\ell(P)$ such that
\begin{align}\label{Vertex-presentation}
v=\left(\bigcap_{i=1}^{d-m}H_i^{'}\right)\cap\left(\bigcap_{j=1}^m H_j^{''}\right),
\end{align}
where $m={\rm dim}\left(\bigcap_{i=1}^{d-m}H_i^{'}\right)$, each $H_i^{'}$ is a facet hyperplane of $\mathcal{O}(P_\ell^{'})$ and each $H_j^{''}$ is a facet hyperplane of $\mathcal{C}(P_{\ell}^{''})$ which corresponds to a chain $C_j$  of length $\geq 2$ in $P_{\ell}^{''}$.  By \cite[Theorem 2.1]{Stanley1986}, there is a set partition $\pi=\{B_1,B_2,\ldots,B_{m+1}\}$ of $[d]$ such that $B_1,B_2,\ldots,B_m$ are connected as subposets of $P_\ell^{'}$,
$B_{m+1}=\{i\in [d]:  a_i=0\ {\rm or}\ 1\}$ and
\begin{align*}
\bigcap_{i=1}^{d-m}H_i^{'}
=&\{(x_1,x_2,\ldots,x_d)\
|\ x_i=x_j\ {\rm\ if} \{i,j\}\subseteq B_k\ {\rm for\ some\ } 1\leq k\leq m,\\
&{\rm and}\ x_r=a_{r}\ {\rm if}\ r\in B_{m+1}\}.
\end{align*}
Let $B_{m+1}=\{r_1,r_2,\ldots,r_s\}$ and for $1\leq k\leq m$, let $b_k$ denote the same values of all $a_i's, i\in B_k$. Then it suffices to show that each $b_k$ is an integer. Keeping in mind the assumption that the Hasse diagram of $P$ is acyclic, we find that $|C_i\cap B_j|\leq 1$ for $1\leq i, j\leq m$. For $1\leq i,j\leq m$, let
\begin{align}\label{Cij}
c_{ij}=\left\{
         \begin{array}{ll}
           1, & \hbox{if $|C_i\cap B_j|=1$;} \\[5pt]
           0, & \hbox{otherwise.}
         \end{array}
       \right.
\end{align}
and for $1\leq i\leq m, 1\leq j\leq s$, let
\begin{align}\label{Dij}
d_{i,m+j}=\left\{
         \begin{array}{ll}
           1, & \hbox{if $r_j\in C_i$;} \\[5pt]
           0, & \hbox{otherwise.}
         \end{array}
       \right.
\end{align}
By \eqref{Vertex-presentation},  $(b_1,b_2,\ldots,b_m,a_{r_1},a_{r_2},\ldots,a_{r_s})$ must be the unique solution of the following linear system:
\begin{align}\label{linear-equations}
\left\{
  \begin{array}{ll}
    \sum_{j=1}^m c_{ij}y_j+\sum_{j=m+1}^{m+s}d_{ij}y_j=1, & 1\leq i\leq m\\[5pt]
    y_{m+1}=a_{r_1}, & \\[5pt]
    y_{m+2}=a_{r_2}, & \\[5pt]
\ \ \ \ \vdots &\\[5pt]
    y_{m+s}=a_{r_s}, & \\[5pt]
  \end{array}
\right.
\end{align}
Now it suffices to show that the determinant of the coefficient matrix
\begin{align}\label{MatrixA}
A=\left(
    \begin{array}{cccccc}
      c_{11} & \cdots & c_{1m} & d_{1,m+1} & \cdots & d_{1,m+s} \\
       &  & \vdots &  &  &  \\
      c_{m1} & \cdots & c_{mm} & d_{m,m+1} & \cdots & d_{m,m+s} \\
      0 & \cdots & 0 & 1 & \cdots & 0 \\
       &  &  & \cdots &  &  \\
      0 & \cdots & 0 & 0 & \cdots & 1 \\
    \end{array}
  \right)
\end{align}
is equal to $1$ or $-1$. Now construct a bipartite graph $G$ with
vertex set
$$\{B_1,B_2,\ldots,B_m,C_1,C_2,\ldots,C_m\}.
$$
and edge set
$$\{\{B_i,C_j\}\ |\ 1\leq i,j\leq m, |B_i\cap C_j|=1\}
$$
Let
\begin{align*}
C=\left(
  \begin{array}{ccc}
    c_{11} & \cdots & c_{1m} \\
     & \vdots &  \\
    c_{m1} & \vdots & c_{mm} \\
  \end{array}
\right)
\end{align*}
Then we have
\begin{align}\label{Det-expansion}
{\rm det}(C)=\sum_{\sigma\in\mathfrak{S_m}}{\rm sign}(\sigma)c_{1\sigma_1}\cdots c_{m\sigma_m}.
\end{align}
Clearly, each nonzero term in \eqref{Det-expansion} corresponds to a perfect matching in the graph $G$. Since the Hasse diagram of $P$ is acyclic, the graph
$G$ must be an acyclic bipartite graph, which means that there is at most one perfect matching in $G$. So we have ${\rm det(C)}=0,1$ or $-1$. Note that the linear equations \eqref{linear-equations} have unique solution $(b_1,b_2,\ldots,b_m,a_{r_1},\ldots,a_{r_s})$. Then we find that ${\rm det}(C)=\pm 1$. It follows that each $b_i$ is an integer. So the vertex $v$ of $\mathcal{OC}_\ell(P)$ is integral. \qed

For general finite poset $P$ with $|E(P)|\geq 2$, the following theorem indicates that there exists at least one nontrivial integral edge partition.

\begin{thm}\label{Existence}
Suppose that $P$ is a finite poset. Let ${\rm Min}(P)$ denote the set of all minimal elements in $P$. For $S\subseteq {\rm Min}(P)$, let $E_S(P)$ denote the set of all edges in
$E(P)$ which are incident to some elements in $S$. Then the edge partition
    $$\ell=(E(P)\setminus E_S(P),E_S(P))$$
 is integral.
\end{thm}

\pf Suppose that $v$ is a vertex of $\mathcal{OC}_\ell(P)$. Then $v$ can be represented as intersection of $d$ independent facet hyperplanes, as in \eqref{Vertex-presentation}. Keeping the notation in the proof of Theorem \ref{THM1}, we can deduce that $|C_i|=2$ and  $|B_i\cap C_j|\leq 1$ for $1\leq i,j\leq m$. So we can construct in the same way two matrices $A$ and $C$ as those in the proof of Theorem \ref{THM1}.  Then, we can construct a graph $G$ with vertex set $\{B_1,B_2,\ldots,B_m, r_1, r_2,\ldots, r_s\}$ and edge set determined by $C_1, C_2,\ldots, C_m.$ More precisely, $\{B_i,B_j\}$ is an edge of $G$ if and only if there exists  $1\leq k\leq m$ such that $C_k=\{i', j'\}$ for some $i'\in B_i, j'\in B_j$, and $\{B_i, r_j\}$ is an edge of $G$ if and only if there exists  $1\leq k\leq m$ such that $C_k=\{r_j, i'\}$ for some $i'\in B_i.$
Obviously, $G$ is a bipartite graph with bipartition $(\mathcal{B}_1,\mathcal{B}_2)$, where
\[
\mathcal{B}_2=\{B_j:1\leq j\leq m, \ B_j=\{k\}\ {\rm for\ some\ } k\in S\}\cup\{r_t:\ 1\leq t\leq s,\ r_t\in S\}.
\]
Moreover, by the construction of the graph $G$, its incidence matrix is
\begin{align*}
 \left(
\begin{array}{cccccc}
     c_{11} & \cdots & c_{1m} & d_{1,m+1} & \cdots & d_{1,m+s} \\
       &  & \vdots &  &  &  \\
      c_{m1} & \cdots & c_{mm} & d_{m,m+1} & \cdots & d_{m,m+s}
      \end{array}
      \right).
\end{align*}
Where $c_{ij}, d_{i,m+j}$ are defined in \eqref{Cij} and in \eqref{Dij} respectively. A well known fact shows that the incidence matrix of any bipartite graph  is totally unimodular. So the
submatrix $C$ has determinant $0, 1$ or $-1$. This completes the proof. \qed

\begin{exa}
By Theorem \ref{THM1}, if the Hasse diagram of $P$ has a cycle, then there exists at least one non-integral edge partition $\ell$.
\begin{itemize}
\item [(1)]
For example, let $P$ denote the poset whose Hasse diagram is a $4$-cycle and let $E_1=\{\{1,2\},\{2,4\},\{3,4\}\}$. Then the edge partition $\ell_1=(E_1,\{1,3\})$ given in Fig. 3(a) is non-integral, since $v=\left(\frac{1}{2},\frac{1}{2},\frac{1}{2}.\frac{1}{2}\right)$ is a vertex of $\mathcal{OC}_{\ell_1}(P)$ given by
 \[
 \left\{
   \begin{array}{ll}
     x_1=x_2=x_4=x_3,&\\
     x_1+x_3=1. &\\
   \end{array}
 \right.
 \]
 Note that the edge partition $\ell_2=(\{1,3\},E_1)$ given in Fig. 3(b) is integral. So we find that the complementary edge partition $\ell^{c}=(cE(P),oE(P))$ of an integral edge partition $\ell=(oE(P),cE(P))$ is not necessarily integral.
\item [(2)] For any poset $P$ whose Hasse diagram is a cycle and any edge partition $\ell$ of $P$, it is not hard to show that all coordinates of each vertex of $\mathcal{OC}_\ell(P)$ are $0$, $1$ or $\frac{1}{2}$.
\end{itemize}
\begin{center}\setlength{\unitlength}{1mm}
\begin{picture}(20,30)
\put(-30,26){$\bullet$}\put(-30,30){$4$}
\put(-38,20){$\bullet$}\put(-42,20){$2$}
\put(-22,20){$\bullet$}\put(-18,20){$3$}
\put(-30,14){$\bullet$}\put(-30,11){$1$}
{\linethickness{0.8mm}\put(-37,21){\line(4,3){8}}
\put(-21,21){\line(-4,3){8}}}
{\linethickness{0.8mm}\put(-29,15){\line(-5,4){8}}}
\put(-29,15){\line(5,4){8}}
\put(-34,5){(a)}
\put(20,26){$\bullet$}\put(20,30){$4$}
\put(12,20){$\bullet$}\put(8,20){$2$}
\put(28,20){$\bullet$}\put(32,20){$3$}
\put(20,14){$\bullet$}\put(20,11){$1$}
\put(13,21){\line(4,3){8}}
\put(29,21){\line(-4,3){8}}
\put(21,15){\line(-5,4){8}}
{\linethickness{0.8mm}\put(21,15){\line(5,4){8}}}
\put(20,5){(b)}
\put(-8,0){Fig. 3}
\end{picture}
\end{center}

\end{exa}

\section{Unimodular equivalence}
In this section, we shall compare the newly constructed order-chain polytopes with some known polytopes. Specifically,  we will focus on integral order-chain polytopes and consider their unimodular equivalence relation with order polytopes or chain polytopes.

We shall use the ideas in the proof of the following theorem due to Hibi and Li \cite{Hibi-Li2012}.

\begin{thm}\cite[Theorem 1.3]{Hibi-Li2012}\label{Hibi-Li}
The order polytope $\mathcal{O}(P)$ and the chain polytope $\mathcal{C}(P)$ of a finite poset $P$ are unimodularly equivalent if and only if the following poset:
\begin{center}\setlength{\unitlength}{1mm}
\begin{picture}(20,40)
\put(-7,15){$\bullet$}\put(-7,10){$1$}
\put(8,15){$\bullet$}\put(8,10){$2$}
\put(-7,30){$\bullet$}\put(-7,33){$4$}
\put(8,30){$\bullet$}\put(8,33){$5$}
\put(0.5,22.5){$\bullet$}\put(6,22.5){$3$}
\put(-6,16){\line(1,1){7}}
\put(1,23){\line(1,1){8}}
\put(9,16){\line(-1,1){7}}
\put(1,24){\line(-1,1){7}}
\put(-5,0){Fig. 4}
\end{picture}
\end{center}
does not appear as a subposet of $P$.
\end{thm}

\begin{defi}
A poset $P$ on $[d]$ is said to be a zigzag poset if its cover relations  are given by
    \[
    1\prec\cdots\prec i_1\succ i_1+1\succ\cdots\succ i_2\prec i_2+1\prec\cdots\prec i_3\succ\cdots\succ i_k\prec i_{k}+1\prec \cdots\prec d
    \] for some $0\leq i_1<i_2<\cdots<i_k\leq d$.
\end{defi}

\begin{thm}\label{Chain}
Suppose that $P$  is a disjoint union of chains. Then for any edge partition $\ell$, the order-chain polytope $\mathcal{OC}_\ell(P)$ is unimodularly equivalent to a chain polytope $\mathcal{C}(Q)$, where $Q$ is a disjoint union of zigzag posets.
\end{thm}
\pf We firstly assume that $P$ is a chain:
\[
1\prec 2\prec 3\prec \cdots \prec d.
\]
and $\ell$ is an edge partition of $P$ given by:
\begin{align*}
{\rm o:}\ \ \ \  &1\prec 2\prec\cdots\prec i_1\\
{\rm c:}\ \ \ \  &i_1\prec i_{1}+1\prec\cdots\prec i_{2}\\
{\rm o:}\ \ \ \  &i_{2}\prec i_{2}+1\prec \cdots\prec i_3\\
\vdots\\
{\rm c:}\ \ \ \   &i_{t-1}\prec i_{t-1}+1\prec \cdots\prec i_t\\
{\rm o:}\ \ \ \   &i_t\prec i_{t}+1\prec \cdots\prec i_{t+1}\\
\vdots\\
{\rm c:}\ \ \ \  &i_{k-1}\prec i_{k-1}+1\prec \cdots\prec i_k=d,
\end{align*}
where $1\leq i_1<i_2<\cdots<i_{k-1}\leq i_k=d$. Then the order-chain polytope $\mathcal{OC}_\ell(P)$ is given by
\begin{align}\label{ORIGINAL}
\left\{
  \begin{array}{ll}
    x_1\geq x_2\geq\cdots\geq x_{i_1}, & \hbox{} \\[7pt]
    x_{i_1}+x_{i_1+1}+\cdots+x_{i_2}\leq 1, & \hbox{} \\[7pt]
    x_{i_2}\geq x_{i_2+1}\geq\cdots\geq x_{i_3}, & \hbox{} \\[7pt]
    \ \ \ \ \ \ \ \ \ \ \ \ \vdots & \hbox{} \\[7pt]
x_{i_{t-1}}+x_{i_{t-1}+1}+\cdots+x_{i_t}\leq 1, &\hbox{}\\[7pt]
x_{i_{t}}\geq x_{i_{t}+1}\geq\cdots\geq x_{i_{t+1}}, &\hbox{}\\[7pt]
\ \ \ \ \ \ \ \ \ \ \ \ \vdots & \hbox{} \\[7pt]
    x_{i_{k-1}}+x_{i_{k-1}+1}+\cdots+ x_{d}\leq 1, & \hbox{}\\[7pt]
    0\leq x_i\leq 1,\ \ \ \ 1\leq i\leq d. &\hbox{}
  \end{array}
\right.
\end{align}
Now define a map $\varphi: \mathbb{R}^{d}\rightarrow \mathbb{R}^{d}$ as follows:
\begin{itemize}
\item[(1)] if $i$ is a maximal element in $P_{\ell}^{'}$, then let $x_i^{'}=x_i$;
\item[(2)] if $i$ is not a maximal element in $P_{\ell}^{'}$, then $\{i,i+1\}$ must be an edge in the Hasse diagram of $P_{\ell}^{'}$. Let $x_i^{'}=x_{i}-x_{i+1}$.
\end{itemize}
Let $\varphi(x_1,x_2,\ldots,x_d)=(x_1^{'},x_{2}^{'},\ldots,x_{d}^{'})$.

Now it is easy to show that $\varphi$ is a unimodular transformation. Moreover, the system \eqref{ORIGINAL} is transformed into:

\begin{align*}
\left\{
  \begin{array}{ll}
    x_1^{'}+x_{2}^{'}+\cdots+x_{i_1}^{'}\leq 1, &\hbox{} \\[8pt]
    x_{i_1}^{'}+x_{i_1+1}^{'}+\cdots+x_{i_2}^{'}+x_{i_2+1}^{'}+\cdots+x_{i_3}^{'}\leq 1, &\hbox{}\\[8pt]
\ \ \ \ \ \ \ \ \ \ \ \ \vdots & \hbox{} \\[8pt]
x_{i_{t-1}}^{'}+x_{i_{t-1}+1}^{'}+\cdots+x_{i_t}^{'}+x_{i_t+1}^{'}+\cdots+x_{i_{t+1}}^{'}\leq 1,\hbox{} \\[8pt]
\ \ \ \ \ \ \ \ \ \ \ \ \vdots & \hbox{} \\[8pt]
    x_{i_{k-1}}^{'}+x_{i_{k-1}+1}^{'}+\cdots+x_{d}^{'}\leq 1, & \hbox{}\\[8pt]
    0\leq x_{i}^{'}\leq 1
  \end{array}
\right.
\end{align*}
Obviously, this system corresponds to the chain polytope $\mathcal{C}(Q)$ for the zigzag poset $Q$:
\begin{align*}
1\prec 2\prec &\cdots\prec i_1\succ i_1+1\succ\cdots\succ i_2\succ i_2+1\succ \cdots\succ i_3\prec\cdots\\
\end{align*}
 or the dual zigzag poset $Q^{*}$:
\begin{align*}
1\succ 2\succ &\cdots\succ i_1\prec i_1+1\prec \cdots\prec i_2\prec i_2+1\prec \cdots\prec i_3\succ\cdots\\
\end{align*}
So we deduce that $\mathcal{OC}_\ell(P)$ is unimodularly equivalent to the chain polytope of  some zigzag poset.

Now we continue to prove the general case that $P$ is a disjoint union of $k$ chains:  $$P=C_1\uplus C_2\uplus\cdots\uplus C_k.$$ Since
\[
\mathcal{O}(P\uplus Q)=\mathcal{O}(P)\times\mathcal{O}(Q)\ {\rm and\ }
\mathcal{C}(P\uplus Q)=\mathcal{C}(P)\times\mathcal{C}(Q),
\]we have
\begin{align}
\nonumber\mathcal{OC}_\ell(P\uplus Q)&=\mathcal{O}((P\uplus Q)_\ell^{'})\cap\mathcal{C}((P\uplus Q)_\ell^{''})\\
\nonumber&=\mathcal{O}(P_\ell^{'}\uplus Q_\ell^{'})\cap\mathcal{C}(P_\ell^{''}\uplus Q_\ell^{''})\\
\nonumber&=\left[\mathcal{O}(P_\ell^{'})\times\mathcal{O}(Q_\ell^{'})\right]\cap
\left[\mathcal{C}(P_\ell^{''})\times\mathcal{C}(Q_\ell^{''})\right]\\
\nonumber&=\left[\mathcal{O}(P_\ell^{'})\cap\mathcal{C}(P_\ell^{''})\right]\times
\left[\mathcal{O}(Q_\ell^{'})\cap\mathcal{C}(Q_\ell^{''})\right]\\
&= \mathcal{OC}_\ell(P)\times\mathcal{OC}_\ell(Q).\label{WITHOUT GENERALITY}
\end{align}
Hence we conclude that
\begin{align*}
\mathcal{OC}_\ell(C_1\uplus\cdots\uplus C_k)&=\mathcal{OC}_\ell(C_1)\times\cdots\times\mathcal{OC}_\ell(C_k)\\
&\overset{\varphi_1\times \cdots\times\varphi_k}{\cong}\mathcal{C}(Q_1)\times\cdots\times\mathcal{C}(Q_k)\\
&=\mathcal{C}(Q_1\uplus\cdots\uplus Q_k),
\end{align*}
where $Q_i$ are zigzag posets. \qed

Similarly, we can modify the proof of Theorem \ref{Chain} slightly to get the following result:

\begin{thm}\label{Zigzag}
Suppose that $P$ is a finite zigzag poset. Then for any edge partition $\ell$, the order-chain polytope $\mathcal{OC}_\ell(P)$ is unimodularly equivalent to a chain polytope $\mathcal{C}(Q)$ for some zigzag poset $Q$.
\end{thm}
 \pf Suppose that $P$ is a zigzag poset on $[d]$ and $\ell$ is an edge partition of $P$. Define a map
$\varphi: \mathbb{R}^{d}\rightarrow\mathbb{R}^{d}$ as follows:
\begin{itemize}
\item[(1)] if $i$ is covered by at most one element in $P_{\ell}^{'}$, let
\begin{align*}
x_{i}^{'}=\left\{
            \begin{array}{ll}
              x_i, & \hbox{if $i$ is a maximal element in\ ;} P_\ell^{'} \\
              x_i-x_j, & \hbox{if $i$ is covered by $j$ in $P_{\ell}^{'}$($j=i-1$ or $i+1$).}
            \end{array}
          \right.
\end{align*}
\item[(2)] if $i$ is covered by both $i-1$ and $i+1$ in $P_{\ell}^{'}$, let
\[
x_i^{'}=1-x_i.
\]
Let $\varphi(x_1,x_2,\ldots,x_d)=(x_1^{'},x_2^{'},\ldots,x_d^{'})$.
\end{itemize}It is not hard to show that $\varphi$ is the desired unimodular transformation.
\qed

The following example shows that not every order-chain polytope $\mathcal{OC}_\ell(P)$ of an acyclic poset $P$ is unimodularly equivalent to some chain polytope.

\begin{exa}\label{NOT-C}
{\em Let $P$ be the poset with an edge partition $\ell$ as follows,
\begin{center}\setlength{\unitlength}{1mm}
\begin{picture}(20,40)
\put(-7,15){$\bullet$}\put(-7,10){$1$}
\put(8,15){$\bullet$}\put(8,10){$2$}
\put(-7,30){$\bullet$}\put(-7,33){$4$}
\put(8,30){$\bullet$}\put(8,33){$5$}
\put(0.5,22.5){$\bullet$}\put(6,22.5){$3$}
{\linethickness{0.8mm}\put(-6,16){\line(1,1){7}}
\put(1,23){\line(1,1){8}}
\put(1,24){\line(-1,1){7}}}
\put(9,16){\line(-1,1){7}}
\put(-5,0){Fig. 5}
\end{picture}
\end{center}
namely, $\ell=(\{\{1,3\},\{3,4\},\{3,5\}\},\{2,3\}).$ Let
\[
\varphi(x_1,x_2,x_3,x_4,x_5)=(x_1,1-x_2,x_3,x_4,x_5).
\]
It is obvious that $\varphi$ is a unimodular transformation and $\varphi(\mathcal{OC}_\ell(P))=\mathcal{O}(P)$. However, by checking all $63$ different
non-isomorphic posets with $5$ elements, we find that $\mathcal{O}(P)$ is not equivalent
to any chain polytope.}
\end{exa}

Furthermore, for any $d\geq 6$, we shall find an integral  order-chain polytope in $\mathbb{R}^d$ which is not unimodularly equivalent to any chain polytope or order polytope. To this end, we need the following lemma.

\begin{lem}\label{nontrivial}
{\rm (1)} None of the chain polytopes of finite posets on $[d]$ possesses
$d + 4$ vertices and $d + 7$ facets.

{\rm (2)} None of the order polytopes of finite posets on $[d]$ possesses
$d + 4$ vertices and $d + 7$ facets.
\end{lem}

\pf
(1)  Assume, by contradiction, that $P$ is a finite poset on $[d]$ such that  $\Cc(P)$ has $d + 4$ vertices  and $d + 7$ facets.
Since the vertices of $\Cc(P)$ are those $\rho(A)$ for which $A$ is an antichain
of $P$, we can deduce that $P$ possesses exactly $d + 4$ antichains.
Keeping in mind that $\emptyset, \{1\}, \ldots, \{d\}$ are antichains of $P$, we find that there is no antichain $A$ in $P$ with $|A| \geq 3$. Otherwise,
the number of antichains of $P$ is at least $d + 5$.
It then follows that there are exactly three $2$-element antichains in $P$.
We need to consider the following four cases:
\begin{itemize}
\item[(i)]
Let, say, $\{1,2\}, \{1,3\}, \{1,4\}$ be the $2$-element antichains of $P$.
Then the maximal chains of $P$ are
$P \setminus \{1\}$ and $P \setminus \{2,3,4\}$.\\
\item[(ii)]
Let, say, $\{1,2\}, \{1,3\}, \{2,4\}$ be the $2$-element antichains of $P$.
Then the maximal chains of $P$ are
$P \setminus \{1,2\}$, $P \setminus \{1,4\}$ and $P \setminus \{2,3\}$.\\
\item[(iii)]
Let, say, $\{1,2\}, \{1,3\}, \{4,5\}$ be the $2$-element antichains of $P$.
Then the maximal chains of $P$ are
$P \setminus \{1,4\}$, $P \setminus \{1,5\}$, $P \setminus \{2,3,4\}$
and $P \setminus \{2,3,5\}$.\\
\item[(iv)]
Let, say, $\{1,2\}, \{3,4\}, \{5,6\}$ be the $2$-element antichains of $P$.
It can be shown easily that $P$ possesses exactly eight maximal chains.
\end{itemize}
Recall that the number of facets of $\Cc(P)$ is equal to
$d + c$, where $c$ is the number of maximal chains of $P$,
it follows from the assumption that there are exactly $7$ maximal chains in $P$, which is a contradiction. As a result, none of the chain polytopes $\Cc(P)$ of a finite poset $P$ on $[d]$ with $d + 4$ vertices can possess $d + 7$ facets, as desired.

\medskip

(2)  Let $P$ be a finite poset on $[d]$ and suppose that the number of vertices
of $\Oc(P)$ is $d + 4$ and the number of facets of $\Oc(P)$ is $d + 7$.
Since the number of vertices of $\Oc(P)$ and that of $\Cc(P)$ coincide,
it follows from the proof of (a) that
there is no antichain $A$ in $P$ with $|A| \geq 3$
and that $P$ includes exactly three $2$-element antichains.
On the other hand, it is known \cite[Corollary 1.2]{Hibi-Li2012} that
the number of facets of $\Oc(P)$ is less than or equal to that of $\Cc(P)$.
Hence the number of maximal chains of $P$ is at least $7$.
Thus, by using the argument in the proof of (a), we can assume that the antichains of $P$ are  $\{1,2\},\{3,4\}$ and $\{5,6\}$.
Then, it is easy to prove that the number of edges of
$\hat{P} = P \cup \{\hat{0}, \hat{1}\}$ is at most $d + 6$,
where $\hat{0} \not\in P$ is the unique minimal element of $\hat{P}$
and $\hat{1} \not\in P$ is the unique maximal element of $\hat{P}$. So we deduce that
 the number of facets of $\Oc(P)$ is at most $d+6$, a contradiction with the assumption.
\qed

We remark that, by modifying the argument of the statement (1) in Lemma 2.6,  we can prove directly that the order polytope of Example \ref{NOT-C}
cannot be unimodularly equivalent to any chain polytope.

\begin{exa}
{\em
Let $P$ be the finite poset as follows (Fig.6):

 \begin{picture}(400,150)(10,50)
 \put(170,170){\circle*{5}}
 \put(170,130){\circle*{5}}
 \put(210,170){\circle*{5}}
 \put(210,130){\circle*{5}}
 \put(210,90){\circle*{5}}
 \put(170,90){\circle*{5}}
 \put(215,168){$6$}
 \put(155,168){$5$}
 \put(215,128){$4$}
 \put(155,128){$3$}
 \put(215,88){$2$}
 \put(155,88){$1$}
 {\linethickness{0.8mm}\put(170,170){\line(0,-1){40}}}
 \put(170,130){\line(0,-1){40}}
 \put(210,170){\line(0,-1){40}}
  \put(210,130){\line(0,-1){40}}
 \put(170,170){\line(1,-1){40}}
 {\linethickness{0.8mm}\put(210,170){\line(-1,-1){40}}}
 \put(170,90){\line(1,1){40}}
 \put(210,90){\line(-1,1){40}}
 \put(170,60){Fig. 6}
 \end{picture}

\noindent
Let $\ell$ be the edge partition with
$oE(P)=\left\{\{3,5\},\{3,6\}\right\}$ and $cE(P)=E(P) \setminus oE(P)$.
Then it is easy to verify that  $\Oc\Cc_{\ell}(P)$ is an integral polytope
with $10$ vertices and $13$ facets.
So it follows from  Lemma \ref{nontrivial} that the integral order-chain polytope
$\Oc\Cc_{\ell}(P)$ cannot be unimodularly equivalent to any order polytope
or any chain polytope.

In fact, for any $d>6$, let $P$ be the following poset and let $\ell$ be the edge partition with $$oE(P)=\left\{\{3,5\},\{3,6\},\{5,7\},\{6,7\},\{7,8\},\cdots,\{d-1,d\}\right\}.$$ It is easy to see that the order-chain polytope $\mathcal{OC}_l(P)$ has $d+4$ vertices and $d+7$ facets. Therefore $\Oc\Cc_{\ell}(P)$ cannot be unimodularly equivalent to any order polytope or any chain polytope.

\begin{picture}(400,210)(10,50)
 \put(170,150){\circle*{5}}
 \put(170,110){\circle*{5}}
 \put(210,150){\circle*{5}}
 \put(210,110){\circle*{5}}
 \put(210,70){\circle*{5}}
 \put(170,70){\circle*{5}}
 \put(215,148){$6$}
 \put(155,148){$5$}
 \put(215,108){$4$}
 \put(155,108){$3$}
 \put(215,68){$2$}
 \put(155,68){$1$}
{\linethickness{0.8mm}\put(170,150){\line(0,-1){40}}}
 \put(170,110){\line(0,-1){40}}
 \put(210,150){\line(0,-1){40}}
  \put(210,110){\line(0,-1){40}}
 \put(170,150){\line(1,-1){40}}
 {\linethickness{0.8mm}\put(210,150){\line(-1,-1){40}}}
 \put(170,70){\line(1,1){40}}
 \put(210,70){\line(-1,1){40}}
 \put(190,180){\circle*{5}}\put(198,175){$7$}
  \put(190,210){\circle*{5}}\put(198,205){$8$}
  {\linethickness{0.8mm}\put(190,180){\line(-2,-3){20}}}
  {\linethickness{0.8mm}\put(190,180){\line(2,-3){20}}}
  {\linethickness{0.8mm}\put(190,180){\line(0,1){30}}}
  {\linethickness{0.8mm}\put(190,210){\line(0,1){10}}}
  \put(188,225){$\vdots$}
    {\linethickness{0.8mm}\put(190,240){\line(0,1){10}}}
     \put(190,250){\circle*{5}}\put(198,245){$d$}
    \put(175,40){Fig.7}
 \end{picture}
}
\end{exa}

Recall that Example $\ref{NOT-C}$ shows that there is an order polytope  which  is not unimodularly equivalent to any chain polytope. To conclude this section, we will prove that, for each $d\geq 9$, there exists a finite poset $P$ on $[d]$ for which the chain
polytope $\mathcal{C}(P)$ cannot be unimodularly equivalent to any order polytope.

Given  a finite poset $P$ on $[d]$,
let $m_\star(P)$ (resp. $m^\star(P)$) denote  the number of minimal
(reps. maximal) elements of $P$ and $c(P)$ denote the number of maximal chains of $P$. For a $d$-dimensional polytope $\mathcal{P}$, denote by $f_{d-1}(\mathcal{P})$ the number of facets of $\mathcal{P}$. Then we have
\[
f_{d-1}(\Oc(P)) = m_\star(P)+m^\star(P)+|E(P)|
\]
and
\[
f_{d-1}(\Cc(P)) = d+c(P),
\]
where $E(P)$ denotes the set of edges
of the Hasse diagram of $P$.
To present our results, we firstly  discuss upper bounds for $f_{d-1}(\Oc(P))$ and $f_{d-1}(\Cc(P))$.
By \cite[Theorem 2.1]{Hibi-Li2012},  if $d \leq 4$, then
$\Oc(P)$ and $\Cc(P)$ are unimodularly equivalent and
$f_{d-1}(\Oc(P)) = f_{d-1}(\Cc(P)) \leq 2d$.
Moreover, for each $1 \leq d \leq 4$, there exists a finite poset $P$ on $[d]$
with $f_{d-1}(\Oc(P)) = f_{d-1}(\Cc(P)) = 2d$.

\begin{lem}
	\label{facet}
Let $d \geq 5$ and $P$ be a finite poset on $[d]$.  Then
\begin{align}\label{UB-O}
f_{d-1}(\mathcal{O}(P))  \leq  \left\lfloor \cfrac{d+1}{2}
\right\rfloor \left( d-\left\lfloor \cfrac{d+1}{2} \right\rfloor \right)+d
\end{align}
and
\begin{align}\label{UB-C}
f_{d-1}(\mathcal{C}(P)) \leq \left\{
\begin{aligned}
&3^k+d, \, \, &d&=3k,\\
&4 \cdot 3^{k-1}+d, \, \, &d&=3k+1,\\
&2 \cdot 3^k+d, \, \, &d&=3k+2.
\end{aligned}
\right.
\end{align}
Furthermore, both  upper bounds for $f_{d-1}(\mathcal{O}(P))$ and $f_{d-1}(\mathcal{C}(P))$  are tight.
\end{lem}
\pf
{\bf (order polytope)}
Let $d = 4$.  Since the right-hand side of
\eqref{UB-O} is equal to $2d \, (= 8)$, the inequality \eqref{UB-O} also holds for $d = 4$. Let $d \geq 5$
and $P$ be a finite poset on $[d]$.
We will prove \eqref{UB-O} by induction on $d$.
Suppose that $1$ is a minimal element of $P$ and let $a$
be the number of elements in $P$ which cover $1$.
	
If $a=0$,  then $\mathcal{O}(P)=\mathcal{O}(P \setminus \{1\})\times [0,1]$
and so
	  \begin{displaymath}
	  \begin{aligned}
	  f_{d-1}(\mathcal{O}(P))&=f_{d-2}(\mathcal{O}(P \setminus \{1\}))+2\\[5pt]
	  &\leq \left\lfloor \cfrac{\,d\,}{2}
	  \right\rfloor\left(d-1-\left\lfloor \cfrac{\,d\,}{2} \right\rfloor\right)+d-1+2\\[5pt]
	  &\leq \left\lfloor \cfrac{d+1}{2}
	  \right\rfloor\left(d-\left\lfloor \cfrac{d+1}{2} \right\rfloor\right)+d.
	  \end{aligned}
	   \end{displaymath}
	
If $1 \leq a \leq \lfloor d / 2 \rfloor$, then  from the facts that
	  $|E(P \setminus \{1\})|=|E(P)|-a$, $m_\star(P\setminus \{1\}) \geq m_\star(P)-1$
	   and $m^\star(P\setminus \{1\}) =m^\star(P)$,
	  we have
	   \begin{displaymath}
	   \begin{aligned}
	   f_{d-1}(\mathcal{O}(P))&=m^\star(P)+m_\star(P)+|E(P)|\\[5pt]
	   &\leq m^\star(P\setminus \{1\})+m_\star(P\setminus \{1\})+1+|E(P\setminus \{1\})| +a \\[5pt]
	   &\leq \left\lfloor \cfrac{d}{2} \right\rfloor\left(d-1-\left\lfloor \cfrac{d}{2}
	   \right\rfloor\right)+(d-1)+\left\lfloor \cfrac{d}{2} \right\rfloor+1\\[5pt]
	   &\leq \left\lfloor \cfrac{d+1}{2} \right\rfloor\left
	   (d-\left\lfloor \cfrac{d+1}{2} \right\rfloor\right)+d.
	   \end{aligned}
	   \end{displaymath}
	
Now we consider the case $\lfloor d/2 \rfloor+1 \leq a \leq d-1$.
Let, say, $2$ be an element of $P$ which covers $1$.
Since the set of the elements of $P$ which cover $1$ is an antichain of $P$,
it follows that $|E(P\setminus\{2\})| \geq |E(P)|-(d-a)$, $m_\star(P\setminus \{2\}) \geq m_\star(P)$
and $m^\star(P\setminus \{2\})  \geq m^\star(P)-1$.  Hence
	      \begin{displaymath}
	      \begin{aligned}
	      f_{d-1}(\mathcal{O}(P))&=m^\star(P)+m_\star(P)+|E(P)|\\[5pt]
	      &\leq m^\star(P\setminus \{2\})+1+m_\star(P\setminus \{2\})+|E(P\setminus \{2\})| +(d-a) \\[5pt]
	      &\leq \left\lfloor \cfrac{\,d\,}{2}
	      \right\rfloor\left(d-1-\left\lfloor \cfrac{\,d\,}{2} \right\rfloor\right)
	      +(d-1)+\left(d-\left\lfloor \cfrac{\,d\,}{2} \right\rfloor-1\right)+1\\[5pt]
	      &\leq \left\lfloor \cfrac{d+1}{2} \right\rfloor\left(d-\left\lfloor \cfrac{d+1}{2}
	      \right\rfloor\right)+d.
	      \end{aligned}
	      \end{displaymath}
	
Therefore, the inequality \eqref{UB-O} holds. We proceed to show that this upper bound for $f_{d-1}(\mathcal{O}(P))$ is tight.  In fact, let $P$ be the finite poset $P$ on $[d]$ with
\[
E(P)=\{ \, \{i,j\} \in [d] \times [d] \, : \,
1 \leq i \leq \left\lfloor \cfrac{d+1}{2} \right\rfloor,
\left\lfloor \cfrac{d+1}{2} \right\rfloor+1 \leq j \leq d \, \}.
\]
Clearly, we have
\[
f_{d-1}(\mathcal{O}(P))=\left\lfloor \cfrac{d+1}{2}
\right\rfloor \left( d-\left\lfloor \cfrac{d+1}{2} \right\rfloor \right)+d.
\]

\medskip

\noindent
{\bf (chain polytope)}
Let $d \geq 5$.  Let $P_{1}$ be a finite poset on $[d]$ and $M_{1}$ the set of
minimal elements of $P_{1}$.
If $P_{1}$ is an antichain, then $f_{d-1}(\mathcal{C}(P_{1})) = 2d$.
Suppose that $P_1$ is not an antichain.
Let $P_2=P_1 \setminus M_1$ and $M_2$ be the set of minimal elements of $P_2$.
In general, if $P_{i}$ is not an antichain and $M_{i}$ is the set of minimal element
of $P_{i}$, then we set $P_{i+1} = P_{i} \setminus M_{i}$.
By continuing this construction, we can get an integer $r \geq 1$ such that
each of the $P_{1}, \ldots, P_{r-1}$ is not an antichain and that $P_{r}$
is an antichain.  Let $P$ be the finite poset on $[d]$ such that
$i_1 \prec i_2 \prec \cdots \prec i_r$ if $i_{j} \in M_{j}$ for $1 \leq j \leq r$.
One has $c(P_1) \leq c(P)=|M_1|\cdots|M_r|$. For any integer $d\geq 5$, let
 \[
M(d)={\rm max}\{\Pi_{i=1}^{r}m_i:\ 1\leq r\leq d,\  m_1+m_2+\cdots+m_r=d, m_i\in\mathbb{N}^{+}\}.
\] Then the desired inequalities \eqref{UB-C} follows
immediately from the following claim:
\begin{align}\label{Maximize the product}
M(d)=\left\{
\begin{aligned}
&3^k, \, \, &d&=3k,\\
&4 \cdot 3^{k-1}, \, \, &d&=3k+1,\\
&2 \cdot 3^k, \, \, &d&=3k+2.
\end{aligned}
\right.
\end{align}
So it suffices to prove this claim. Since for any integer $m \geq 4$,
\[
m \leq \left\lfloor \cfrac{m+1}{2} \right\rfloor
\left(m-\left\lfloor \cfrac{m+1}{2} \right\rfloor \right),
\]
we can assume that, to maximize the product $\Pi_{i=1}^{r}m_i$, all parts $m_i\leq 3$.
We can also assume without loss of generality  that  there are at most  two $m_i$s that are equal to $2$, since $2^3<3^2$. Then the claim \eqref{Maximize the product} follows immediately.

Finally, for each $d \geq 5$,
the existence of a finite poset $P$ on $[d]$ for which the equality holds
in \eqref{UB-C} follows easily from the above argument.
\qed

A routine computation shows that,
for each $1 \leq d \leq 8$,
the right-hand side of \eqref{UB-O} coincides with that of \eqref{UB-C} and that,
for each $d \geq 9$, the right-hand side of \eqref{UB-O} is strictly less than
that of \eqref{UB-C}.  Hence

\begin{cor}
For each $d \geq 9$, there exists a finite poset $P$ on $[d]$ for which
the chain polytope $\Cc(P)$ cannot be unimodularly equivalent to any order polytope.
\end{cor}

\section{Volumes of $\mathcal{OC}_\ell(P)$}

Given a poset $P$ on $[n]$, Corollary 4.2 in \cite{Stanley1986} shows that the volumes of $\mathcal{O}(P)$ and $\mathcal{C}(P)$ are given by
\[
V(\mathcal{O}(P))=V(\mathcal{C}(P))=\frac{e(P)}{n!},
\]
where $e(P)$ is the number of linear extensions of $P$. (Recall that a linear extension of $P$  is a permutation $\pi=\pi_1\pi_2\cdots\pi_n$ of $[n]$ such that $\pi^{-1}(i)<\pi^{-1}(j)$ if $i\prec j$ in $P$.)

For order-chain polytopes, different edge partitions usually give rise to polytopes with different volumes. For example, let $P$ be the poset as follows.
  \newline
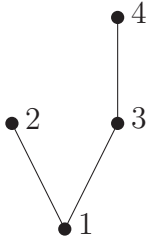
  \begin{picture}(400,150)(10,50)
  \put(210,150){\circle*{5}}
  \put(170,110){\circle*{5}}
  \put(210,110){\circle*{5}}
  \put(190,70){\circle*{5}}
       \put(215,108){$3$}
       \put(215,148){$4$}
       \put(195,68){$1$}
       \put(175,108){$2$}
  \put(210,110){\line(0,1){40}}
    \put(190,70){\line(-1,2){20}}
     \put(190,70){\line(1,2){20}}

  \end{picture}\\
It is easy to see that $$V(\Oc(P))=V(\Cc(P))=\frac{3}{4!}.$$ Let \begin{align*}\ell=(\{1,2\},\{\{1,3\},\{3,4\}\}),\ \  \ \ \ \ \ell'=(\{\{1,2\},\{1,3\}\},\{3,4\}),
\end{align*}
 then  we have
 \begin{align*}
 V(\Oc\Cc_\ell(P))=\frac{1}{4!}\ \ \ \ {\rm and}\ \ \ \ \
   V(\Oc\Cc_{\ell'}(P))=\frac{5}{4!}.
   \end{align*}
  Hence one has the following inequality:
  $$ V(\Oc\Cc_\ell(P)) < V(\Oc(P))=V(\Cc(P))<V(\Oc\Cc_{\ell'}(P)).$$
Then a natural question is to ask which edge partition $\ell$  gives rise to an order-chain polytope with maximum volume. It seems very difficult to solve this problem in general case. In this section, we consider the special case when $P$ is a chain $P$ on $[n]$. We transform it to a  problem of maximizing descent statistics over certain family of subsets. For references on this topic,  we refer the reader to \cite{Ehrenborg-Mahajan1998} and \cite{Sagan-Yeh-Ziegler2013}.

Let $P$ be a chain on $[n]$. By the proof of Theorem \ref{Chain}, for an edge partition $\ell$ of $P$, the order-chain polytope $\mathcal{OC}_\ell(P)$ is unimodularly equivalent to a chain polytope $\mathcal{C}(P_1)$, where $P_1$ is a zigzag poset such that all maximal chains,  except the first one (containing $1$) and the last one (containing $n$), consist of at least three elements. So we have
\[
V(\mathcal{OC}_{\ell}(P))=V(\mathcal{C}(P_1))=\frac{e(P_1)}{n!}.
\]
Conversely, for such a zigzag poset $P_1$, it is easy to find an edge partition $\ell$ of $P$ such that $\mathcal{OC}_\ell(P)$ is unimodularly equivalent to $\mathcal{C}(P_1).$ Denote by $\mathcal{Z}(n)$ the set of such zigzag posets $P_1$ on $[n]$. Thus, to compute the maximum volume over all order-chain polytopes of the chain $P$, it suffices to compute the maximum number of linear extensions for all  zigzag posets $P_1\in\mathcal{Z}(n)$. Next we shall represent this problem as a problem of maximizing
 descent statistic over a certain class of subsets. To this end, we recall some notions and basic facts. Given a permutation $\pi=\pi_1\pi_2\cdots\pi_n$, let  ${\it Des}(\pi)$ denote its descent set $\{i\in [n-1]: \pi_i>\pi_{i+1}\}$. For $S\subseteq [n-1]$, define the descent statistic $\beta(S)$ to be the number of permutations of $[n]$ with descent set $S$. Note that there is an obvious bijection between zigzag posets on $[n]$ and subsets of $[n-1]$ given by
$$S: P\mapsto \{j\in [n-1]: j\succ j+1\}.$$
Moreover,  a permutation $\pi=\pi_1\pi_2\cdots\pi_n$ of $[n]$ is a linear extension of $P$ if and only if $\rm{Des}(\pi^{-1})=S(P)$. Let $\mathcal{F}(n)=S(\mathcal{Z}(n))$. Then we can transform the problem of maximizing volume of order-chain polytopes of an $n$-chain  to the problem of maximizing the descent statistic $\beta(S)$, where $S$ ranges over $\mathcal{F}(n)$.

Observe that $\beta(S)=\beta(\bar{S})$, where $\bar{S}=[n-1]\setminus S$.  Following \cite{Ehrenborg-Mahajan1998},  we will encode both $S$ and $\bar{S}$ by a list $L=(l_1,l_2,\ldots,l_k)$ of positive integers such that $l_1+l_2+\cdots+l_k=n-1$. Given $S\subseteq [n-1]$, a run of $S$ is a set $R\subseteq [n-1]$ of consecutive integers of maximal cardinality such that $R\subseteq S$ or $R\subseteq \bar{S}$. For example, if $n=10$, then the set $S=\{1,2,5,8,9\}$ has $5$ runs: $\{1,2\},\ \{3,4\},\ \{5\},\ \{6,7\},\ \{8,9\}.$   Suppose that $S$ has $k$ runs $R_1, R_2,\ldots, R_k$ with $|R_i|=l_i$,  let $L(S)=(l_1,l_2,\ldots,l_k)$.

\begin{lem}\label{list}
Suppose that $S\subseteq [n-1]$ and $L(S)=(l_1,l_2,\ldots l_k)$. Then $S\in\mathcal{F}(n)$ if and only if $l_i\geq 2$ for all $2\leq i\leq k-1$.
\end{lem}
\pf The lemma follows immediately from the fact that $\mathcal{Z}(n)$ consists of zigzag posets $P$ such that all maximal chains in $P$,  except the first one (containing $1$) and the last one (containing $n$), contains at least three elements.  \qed

Denote by $F_n$ the $n$th Fibonacci number. By Lemma \ref{list}, it is easy to see that $|\mathcal{F}(n)|=2F_n$ for $n\geq 2$. Based on computer evidences, we conjectured the following results about maximizing descent statistic over $\mathcal{F}(n)$, which in fact \footnote{We thank Joe Gallian and  Mitchell Lee  for bringing \cite[Theorem 6.1]{Ehrenborg-Mahajan1998} to our attention.} is a special case of Theorem \cite[Theorem 6.1]{Ehrenborg-Mahajan1998}.

\begin{prop} Suppose that $n\geq 2$ and $S\subseteq [n-1]$.
\begin{itemize}
\item[(1)] If $n=2m$  and
\begin{align*}
L(S)=(1,\underbrace{2,2,\ldots,2}_{m-1})\ \ {\it or}\ \ L(S)=(\underbrace{2,2,\ldots,2}_{m-1},1),
\end{align*}
then $\beta(T)\leq \beta(S)$ for any $T\in \mathcal{F}(n)$.\\[5pt]
\item[(2)] If $n=2m+1$  and
\begin{align*}
L(S)=(1,\underbrace{2,2,\ldots,2}_{m-1},1),
\end{align*}then
$\beta(T)\leq \beta(S)$ for any $T\in\mathcal{F}(n)$.
\end{itemize}
\end{prop}

Equivalently, by the proof of Theorem \ref{Chain}, we have

\begin{prop}
Let $P$ be a chain on $[n]$. Then the alternating edge partition $\ell=(oE(P),cE(P))$ with
\begin{align*}
oE(P)=\left\{
        \begin{array}{ll}
          \{\{1,2\},\{3,4\},\ldots,\{n-1,n\}\}, & \hbox{ if $n$ is even;} \\[7pt]
          \{\{1,2\},\{3,4\},\ldots,\{n-2,n-1\}\}, & \hbox{otherwise.}
        \end{array}
      \right.
\end{align*}
gives rise to an order-chain polytope $\mathcal{OC}_\ell(P)$ with maximum volume.
\end{prop}

\vspace{0.5cm}
 \noindent{\bf Acknowledgments.}  This work was initiated when the third author and the fourth author were visiting the Math department of M.I.T. These two authors would like to thank Professor Richard Stanley for many helps and the whole Math department of M.I.T for providing a great environment.  This work was supported by the Research Foundation for the Doctoral Program of Higher Education of China (Grant No. 20130182120030) and the China Scholarship Council.

\end{document}